\documentclass[11pt]{article}
\usepackage{graphicx}

\usepackage{amssymb}	
\usepackage{amsmath}   
\usepackage{latexsym}	
\usepackage{longtable}
\usepackage{multirow}
\usepackage{epsfig}

\allowdisplaybreaks

\newtheorem{algorithm}{Algorithm}[section]

\newtheorem{theorem}{Theorem}[section]

\newtheorem{remark}{Remark}[section]

\newcommand{\qed}{\nobreak \ifvmode \relax \else \ifdim\lastskip<1.5em \hskip-\lastskip \hskip1.5em plus0em minus0.5em \fi \nobreak \vrule height0.75em width0.5em depth0.25em\fi} 

\def\A{{\bf A}}
\def\B{{\bf B}}
\def\C{{\bf C}}
\def\D{{\bf D}}

\def\F{{\bf F}}
\def\G{{\bf G}}
\def\H{{\bf H}}
\def\I{{\bf I}}
\def\J{{\bf J}}
\def\K{{\bf K}}

\def\0{{\bf 0}}

\def\Q{{\bf Q}}
\def\R{{\bf R}}
\def\S{{\bf S}}
\def\T{{\bf T}}

\def\V{{\bf V}}

\def\X{{\bf X}}
\def\Y{{\bf Y}}

\def\a{{\bf a}}

\def\d{{\bf d}}
\def\e{{\bf e}}
\def\f{{\bf f}}

\def\h{{\bf h}}

\def\n{{\bf n}}
\def\p{{\bf p}}
\def\q{{\bf q}}
\def\r{{\bf r}}

\def\t{{\bf t}}
\def\u{{\bf u}}

\def\x{{\bf x}}
\def\y{{\bf y}}

\def\Tr{{\rm T}}

\def\Re{{\cal R}_e}

\def\diag{{\rm diag}}

\title{Model Predictive Control in Spacecraft Rendezvous and Soft Docking}
\author {Yaguang Yang\thanks{
Office of Research, NRC, 21 Church Street, Rockville, 20850. Email:
yaguang.yang@verizon.net} 
}
\date{\today}

\begin{document}

\maketitle  

\begin{abstract}
This paper discusses translation and attitude control in 
spacecraft rendezvous and soft docking. The target spacecraft
orbit can be either circular or elliptic. The high fidelity
model for this problem is intrinsically a nonlinear system
but can be viewed as a linear time-varying system (LTV). 
Therefore, a model predictive control (MPC) based design 
is proposed to deal with the time-varying
feature of the problem. A robust pole assignment method is
used in the MPC-based design because of the following merits
and/or considerations: (a) no overshoot of the relative
position and attitude between the target and the chaser to 
achieve soft docking by placing all closed-loop poles
in the negative real axis of the complex plan, which avoids
oscillation of the relative position and attitude, in 
particular, in the final stage, (b) fast on-line 
computation, (c) modeling error tolerance, and (d) disturbance 
rejection. We will discuss these considerations and merits, and 
use a design simulation to demonstrate that the desired 
performance is indeed achieved.
\end{abstract}


\newpage

\section{ Introduction}

Spacecraft rendezvous is an important operation in many space 
missions. There are extensive research in this field and hundreds 
successful rendezvous missions, see, for example, the survey 
paper \cite{lzt14} and references therein. The entire rendezvous 
process can be divided into several phases, including phasing, 
close-range rendezvous, final approaching, and docking. In the 
early phase, the chaser fly to the target with the aid from the 
ground station and translation control is the main concern. For
this purpose, the well-known Hill \cite{hill1878} or Clohessy 
and Wiltshire \cite{cw60} equations are adequate for the control 
system design if the orbit is circular. But in the final 
approaching and docking phase, both translation and attitude 
control may be required. Moreover, we may need to consider the
case that the orbit of the target spacecraft is not circular. 
To achieve this requirement, more complex models, 
for example, those discussed in
\cite{wh96,pk01,kgng07}, should be considered. Although these
models are for more general purpose, they can be easily tailored 
for the use of spacecraft rendezvous and docking control.

The research of spacecraft rendezvous has attracted renewed 
interest in recent years as a result of new development in 
control theory and increased space missions involving rendezvous
and soft docking. Various design methods have been considered
for this control system design problem. For example,
an adaptive output feedback control was proposed for this purpose
in \cite{ssj06}; a multi-objective robust $H_{\infty}$ control method
was investigated in \cite{gys09};  a Lyapunov differential equation 
approach was studied for elliptical orbital rendezvous with constrained 
controls \cite{zld11}; a gain scheduled control of linear 
systems was applied to spacecraft rendezvous problem subject
to actuator saturation \cite{zwld14}; and various control design
methods were considered for 6 degree of freedom (DOF) spacecraft 
proximity operations \cite{sh15,zd12,kng08,slf11,xjlky15}.
All these methods have their merits in solving the challenging 
problem under various conditions, but none of them addressed
a fundamental issue, i.e., to achieve the soft docking.

In this paper, we first carefully examine a newly proposed
model derived in \cite{kgng07}, then we determine the measurable 
variables and controllable inputs in the mission
of the final approaching and docking phase, and make some 
reasonable assumptions that normally hold via engineering
design. We also adopt a reduced quaternion concept proposed in
\cite{yang10} to slightly simplify the model (because of some 
merit discussed in \cite{yang10, yang14}). To make the general
model useful for the control system design, thruster configuration is
considered and modeled in a complete control system model. This complete 
model can be viewed either as a nonlinear model or a linear 
time-varying (LTV) model. We prefer to use the linear time-varying
model because a linear system is easier to handle than a nonlinear
system and the corresponding design
methods are capable to consider the system performance which
is very important as  {\it soft docking may not 
allow overshoot in the relative position in the spacecraft 
rendezvous and docking phase}.

There are two popular methods that deal with linear time-varying control 
system design with the consideration of system performance. The 
first one is gain scheduling \cite{rs00,zwld14} and the second one is 
model predictive control \cite{aw13}. A simple analysis in 
\cite{yang18} shows that the former is the most efficient when 
all time-varying parameters explicitly depends on time; and the 
later is more appropriate when many parameters depend implicitly 
on time. The rendezvous and docking model falls into the second
category. Therefore, we propose a MPC-based method to design the rendezvous and docking control. Although several design methods, 
such as LQR, $\H_{\infty}$, and robust pole assignment, take the 
performance into the design consideration and can be the candidates
for the MPC-based design, only robust pole
assignment method can directly take overshoot into the design
consideration because overshoot is directly related to the 
closed-loop pole positions \cite{dh08}. In addition, robust
pole assignment guarantees that the closed-loop poles are
not sensitive to the parameter changes in the system 
\cite{psnyst14} that is important given that the system is time-varying.  
Moreover, robust pole assignment design minimizes an upper bound 
of $\H_{\infty}$ norm which means that the design is robust to 
the modeling error and reduces the impact of disturbance torques 
on the system output \cite{yang14}. Among many robust pole 
assignment algorithms, we suggest a globally convergent 
algorithm proposed in \cite{tity96} because of its fast on-line computation 
and other merits \cite{psnyst14}. We indicate that the closed-loop
system is exponentially stable and use a design example 
to show the efficiency and effectiveness of the 
proposed method.


The remainder of the paper is organized as follows. Section 2 
summarizes the complete rendezvous model of \cite{kgng07},
simplifies the model for the use for rendezvous and soft docking
control, and provides engineering design considerations
for rendezvous and docking control. Section
3 discusses the MPC-based method for spacecraft control
using robust pole assignment. Section 4 provides a design example
and simulation result. Section 5 is the summary of conclusions
of the paper.

\section{Spacecraft model for rendezvous} 

In this section, we first introduce the model developed by 
Kristiansen et. al. in \cite{kgng07}. We also discuss the 
assumptions required from the application of final approaching and 
docking phase in the rendezvous process and derive a simplified 
model to be used in this paper. 
For the sake of simplicity, we denote a vector by bold 
symbol $\a$ and its magnitude or 2-norm $\| \a \|$  by normal font $a$. 
We make the following assumption throughout the paper.
\vspace{0.05in}
\newline
{\bf Assumption 1:} {\it Chaser and target can exchange 
information such as position,
attitude and rotational rate in real time.}
\vspace{0.03in}

This assumption can be achieved by engineering design.

\subsection{The model for translation dynamics}

\begin{figure}[h!]
\centering
\includegraphics[height=14cm,width=10cm]{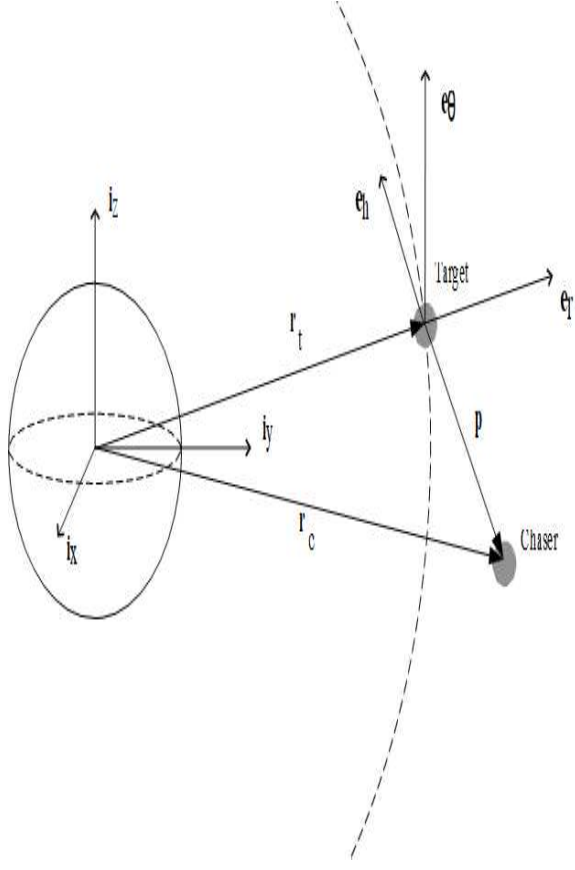}
\caption{Coordinate Frame.}
\label{fig:coordinate}
\end{figure}

As shown in Figure \ref{fig:coordinate}, the inertial frame is 
defined by standard earth-centered inertial (ECI) frame 
${\cal F}_i$ \cite{yang12}. Let $\r_t$ be the vector from
the Earth center to the center of the mass of the target and 
$\r_c$ be the vector from
the Earth center to the center of the mass of the chaser.
Let the angular momentum vector of the target orbit be denoted by 
$\h= \r_t \times \dot{\r}_t$. The orbital frame is the standard 
local vertical local horizontal (LVLH) frame ${\cal F}_t$ 
centered on the target with the basis vectors of the frame 
$\e_r =\r_t/ r_t$, $\e_h = \h/h$, and $\e_{\theta}=\e_h \times \e_r$.
The relative position vector between a target and a chaser is defined
by 
\begin{equation}
\p=\r_c - \r_t = x \e_r + y \e_{\theta} + z \e_h.
\label{relativeP}
\end{equation}
$\p$ is available in real time if GPS is installed in both 
spacecraft and Assumption 1 holds.
The body frames of the target and the chaser, ${\cal F}_{tb}$  and 
${\cal F}_{cb}$, have their origins at their centers of mass and
their coordinate vectors are their principal axes of the inertia.
Let $\gamma$ be the true anomaly of the target, $m_t$ be the 
mass of the target and $m_c$ be the mass of the chaser, 
$\mu$ be the geocentric gravitational constant. Then, the 
nonlinear {\it relative position dynamics} can be expressed in the
target frame as \cite{kgng07}:
\begin{equation}
m_c \ddot{\p} + \C_t(\dot{\gamma}) \dot{\p} + \D_t(\dot{\gamma}, \ddot{\gamma}, r_c) \p +\n_t(r_c,r_t)
= \f_c + \f_d,
\label{translationDynamics}
\end{equation}
where 
\begin{equation}
\C_t(\dot{\gamma}) =2m_c \dot{\gamma} \left[ \begin{array}{ccc}
0 & -1 & 0 \\
1 & 0 & 0 \\
0 & 0 & 0 
\end{array}  \right]
\label{cMatrix}
\end{equation}
\begin{equation}
\D_t(\dot{\gamma}, \ddot{\gamma}, r_c)
 = m_c \left[ \begin{array}{ccc}
\frac{\mu}{r_c^3} -\dot{\gamma} & -\ddot{\gamma} & 0 \\
\ddot{\gamma} & \frac{\mu}{r_c^3} -\dot{\gamma}  & 0 \\
0 & 0 & \frac{\mu}{r_c^3} 
\end{array}  \right]
\label{dMatrix}
\end{equation} 
\begin{equation}
\n_t(r_c,r_t)=m_c \mu \left[ r_t/r_c^3 - 1/r_t^2, 0, 0 \right]^{\Tr},
\label{nVector}
\end{equation}
$\f_c$ is the control force vector, and $\f_d$ is the 
disturbance force vector, both are applied in chaser's body 
frame. It is worthwhile to note that 
\begin{equation}
\n_t(r_c,r_t) \big|_{r_c=r_t}=\0.
\label{ntZero}
\end{equation} 
Also, according to Assumption 1, $\C_t(\dot{\gamma})$,
$\D_t(\dot{\gamma}, \ddot{\gamma}, r_c)$
and $\n_t(r_c,r_t)$ are known but are time-varying.

\subsection{The model for attitude dynamics}

Let the unit quaternion $\bar{\q} = \left[ q_0, \q^{\Tr} \right]^{\Tr}$
be the relative attitude of the target and chaser, where
\begin{equation}
\q^{\Tr} = [q_1, q_2, q_3].
\end{equation}
The inverse of the quaternion is defined in \cite{yang12} as
$\bar{\q}^{-1}  = \left[ q_0, - \q^{\Tr} \right]^{\Tr}$.
Let $\bar{\q}_{i,cb}=[q_{c0}, q_{c1}, q_{c2}, q_{c3}]$ be the 
relative quaternion of from chaser's body frame to the inertial 
frame, and $\bar{\q}_{i,tb}=[q_{t0}, q_{t1}, q_{t2}, q_{t3}]$ 
be the the relative quaternion of from 
target's body frame to the inertial frame. Notice that 
$\bar{\q}_{i,cb}$ is measurable from the chaser and 
$\bar{\q}_{i,tb}$ is measurable from the  target.
Using the Assumption 1, equation (7.2.14) of \cite{sidi97},
and equation (45b) of \cite{yang12}, we have
\begin{equation}
\bar{\q} = \bar{\q}_{i,cb}^{-1} \bar{\q}_{i,tb}
=\left[ \begin{array}{cccc}
q_{t0} & -q_{t1} & -q_{t2} & -q_{t3} \\
q_{t1} & q_{t0} & q_{t3} & -q_{t2} \\
q_{t2} & -q_{t3} & q_{t0} & q_{t1} \\
q_{t3} & q_{t2} & -q_{t1} & q_{t0} 
\end{array} \right]
\left[ \begin{array}{c}
q_{c0}  \\
-q_{c1} \\
-q_{c2} \\
-q_{c3} 
\end{array} \right],
\label{relativeQ}
\end{equation}
which, according to Assumption 1, is measurable.
The relative angular velocity between frames ${\cal F}_{cb}$ and 
${\cal F}_{tb}$ expressed in frame ${\cal F}_{cb}$ is given by
\begin{equation}
\boldsymbol{\omega} = \boldsymbol{\omega}_{i,cb}^{cb}-
\R_{tb}^{cb} \boldsymbol{\omega}_{i,tb}^{tb}
=[\omega_1, \omega_2, \omega_3]^{\Tr},
\label{relativeW}
\end{equation}
where $\boldsymbol{\omega}_{i,cb}^{cb}$ is the angular velocity of the 
chaser body frame relative to the inertial frame expressed in the chaser
body frame, which is measurable from the
chaser; $\boldsymbol{\omega}_{i,tb}^{tb}$ is the angular velocity of the 
target body frame relative to the inertial frame expressed in the target
body frame, which is measurable from the
target; $\R_{tb}^{cb}$ is the rotational matrix from ${\cal F}_{tb}$
to ${\cal F}_{cb}$ which is an equivalent rotation of $\bar{\q}$ 
and can be expressed as \cite[eq. (43)]{yang12}
\begin{equation}
\R_{tb}^{cb} = (q_0^2-\q^{\Tr} \q) \I + 2 \q \q^{\Tr}- 2q_0 \S(\q),
\label{relativeFrame}
\end{equation}
where $\S( \q )= \q {\times}$ is a cross product operator. Using
Assumption 1 again, we conclude that $\boldsymbol{\omega}$ 
is available from measurements. Let $\J_c$ and $\J_t$ be 
the inertia matrices of the chaser and the target, respectively. 
The relative attitude dynamics is given in chaser's frame as
\cite{kgng07,zd12}:
\begin{equation}
\J_c \dot{\boldsymbol{\omega}} + \C_r(\boldsymbol{\omega},\q) 
\boldsymbol{\omega} + \n_r(\boldsymbol{\omega},\q)
= \t_c + \t_d,
\label{relativeDynamics}
\end{equation}
where $\t_c$ and $\t_d$ are control torque and disturbance 
torque respectively, both are applied to the chaser's body frame, 
$\C_r(\boldsymbol{\omega})$ and $\n_r(\boldsymbol{\omega})$ 
are given as follows:
\begin{equation}
\C_r(\boldsymbol{\omega},\q) =\J_c \S(\R_{tb}^{cb} \boldsymbol{\omega}_{i,tb}^{tb})
+\S(\R_{tb}^{cb} \boldsymbol{\omega}_{i,tb}^{tb}) \J_c
-\S(\J_c(\boldsymbol{\omega} +\R_{tb}^{cb} \boldsymbol{\omega}_{i,tb}^{tb})),
\label{Cr}
\end{equation}
\begin{equation}
\n_r(\boldsymbol{\omega},\q)
=\S(\R_{tb}^{cb} \boldsymbol{\omega}_{i,tb}^{tb})
\J_c  \R_{tb}^{cb} \boldsymbol{\omega}_{i,tb}^{tb}
-\J_c  \R_{tb}^{cb}  \J_t^{-1}\S(\boldsymbol{\omega}_{i,tb}^{tb})
\J_t \boldsymbol{\omega}_{i,tb}^{tb}).
\label{nr}
\end{equation}
At the end of the docking phase, the rotation matrix satisfies
$\R_{tb}^{cb} =\I$. 

For attitude dynamics, we use the reduced
quaternion as proposed in \cite{yang10} which is given as follows\footnote{
In the final approaching and docking phase,  the attitude
error between chaser and target is very small, the quaternion 
is far away from the only singular point of the reduced quaternion.}:
\begin{eqnarray} 
\dot{\q} & = & \left[  \begin{array} {c} \dot{q}_1 \\ \dot{q}_2 \\ \dot{q}_3
\end{array} \right]
= \frac{1}{2}  \left[  \begin{array} {ccc} 
\sqrt{1-q_1^2-q_2^2-q_3^2} & -q_3  & q_2 \\
q_3 & \sqrt{1-q_1^2-q_2^2-q_3^2} & -q_1 \\
-q_2  &  q_1  & \sqrt{1-q_1^2-q_2^2-q_3^2} \\
\end{array} \right] 
\left[  \begin{array} {c}  \omega_1 \\ \omega_2 \\ \omega_3
\end{array} \right] 
\nonumber  \\
& = & \frac{1}{2}  \T \boldsymbol{\omega} .
\label{qdot}
\end{eqnarray}

\subsection{A complete model for rendezvous and docking}

Let 
\begin{equation}
{\r}= \dot{\p},
\label{rdot}
\end{equation}
which can be approximated by $\dot{\p} \approx \Delta \p / \Delta t$. 
Now, we can summarize the result by combining equations
(\ref{rdot}), (\ref{translationDynamics}), 
(\ref{relativeDynamics}), and (\ref{qdot}), which yields
\begin{equation}
\dot{\x} = 
\left[ \begin{array}{c}
\dot{\p}  \\
\dot{\r}  \\
\dot{\q}  \\
\dot{\boldsymbol{\omega}}  
\end{array} \right]
=\left[ \begin{array}{cccc}
\0 & \I & \0  & \0  \\
-\frac{1}{m_c} \D_t & -\frac{ 1}{m_c} \C_t &  \0  &  \0   \\
 \0   &  \0   &  \0   &  \frac{1}{2}  \T  \\
 \0    & \0    &  \0    & -\J_c^{-1} \C_r
\end{array} \right]
\left[ \begin{array}{c}
{\p}  \\
{\r}   \\
{\q}   \\
{\boldsymbol{\omega}} 
\end{array} \right]
-\left[ \begin{array}{c}
 \0    \\
\frac{1}{m_c} \n_t   \\
 \0     \\
 \J_c^{-1} \n_r 
\end{array} \right]
+\left[ \begin{array}{c}
 \0    \\
\frac{1}{m_c}  \f_c  \\
 \0     \\
 \J_c^{-1} \t_c 
\end{array} \right].
\label{rendezvousModel}
\end{equation}
Since $\D_t$, $\C_t$, $\T $, $\C_r$, $\n_t$, and $\n_r$ depend 
on $\q$, $\boldsymbol{\omega}$, $r_c$, $r_t$, $\gamma$ which 
are all time-varying, equation (\ref{rendezvousModel}) can be 
treated as a linear time-varying system.

\begin{figure}[hb!]
\centering
\includegraphics[height=9cm,width=7cm]{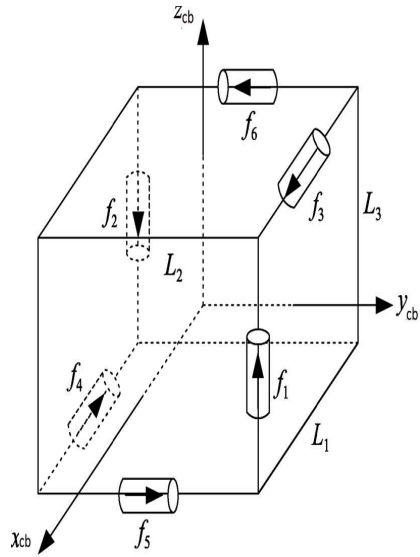}
\caption{Thruster configuration.}
\label{fig:thrusterInstallation}
\end{figure}

For spacecraft rendezvous and docking, the control forces
and torques are provided by thruster systems.
It is well-known that the control force vector and control 
torque vector depend on the thruster configurations and many
configurations are reported in different systems, for example,
\cite{xtf12,crb10,yang14b}. Let $\F_a$ and $\T_a$ be the 
thruster configuration related matrices that define the control 
force vector and control torque vector, i.e.,
\begin{equation}
\f_c = \F_a \f_a, \hspace{0.1in} \t_c = \T_a \f_a,
\label{controlForceV}
\end{equation}
where $\f_a$ is the vector of forces generated by thrusters. 
Substituting (\ref{controlForceV}) into
(\ref{rendezvousModel}), we have 
\begin{eqnarray}
\dot{\x} & = &
\left[ \begin{array}{c}
\dot{\p}  \\
\dot{\r}  \\
\dot{\q}  \\
\dot{\boldsymbol{\omega}}  
\end{array} \right]
=\left[ \begin{array}{cccc}
\0 & \I & \0  & \0  \\
-\frac{1}{m_c} \D_t & -\frac{ 1}{m_c} \C_t &  \0  &  \0   \\
 \0   &  \0   &  \0   &  \frac{1}{2}  \T  \\
 \0    & \0    &  \0    & -\J_c^{-1} \C_r
\end{array} \right]
\left[ \begin{array}{c}
{\p}  \\
{\r}   \\
{\q}   \\
{\boldsymbol{\omega}} 
\end{array} \right]
-\left[ \begin{array}{c}
 \0    \\
\frac{1}{m_c} \n_t   \\
 \0     \\
 \J_c^{-1} \n_r 
\end{array} \right]
+\left[ \begin{array}{c}
 \0    \\
\frac{1}{m_c}  \F_a  \\
 \0     \\
 \J_c^{-1} \T_a
\end{array} \right] \f_a
\nonumber \\
& = & \A(t) \x - \n_d(t) + \B \f_a.
\label{rendezvousModel2}
\end{eqnarray}
Assuming that the chaser's mass change due to fuel consumption 
is negligible, the matrix $\B$ is then a constant. For illustrative
purpose, in the rest of the discussion, we assume that the thrusters have 
the configuration considered in \cite{zd12} which is described
in Figure \ref{fig:thrusterInstallation}. But the same idea
can be used in other thruster configurations. Therefore, we have

\begin{equation}
\F_a  =\left[ \begin{array}{cccccc}
0 & 0 & 1 & -1 & 0 & 0 \\
0 & 0 & 0 & 0 & 1 & -1 \\
1 & -1 & 0 & 0 & 0 & 0 
\end{array} \right],
\label{forceConfig}
\end{equation}
and 
\begin{equation}
\T_a  =\left[ \begin{array}{cccccc}
\frac{L_2}{2} & \frac{L_2}{2} & 0 & 0 & \frac{L_3}{2} & \frac{L_3}{2}  \\
-\frac{L_1}{2} & -\frac{L_1}{2} & \frac{L_3}{2} & \frac{L_3}{2}  & 0 & 0 \\
0 & 0 & -\frac{L_2}{2} & -\frac{L_2}{2} & \frac{L_1}{2} & \frac{L_1}{2} 
\end{array} \right].
\label{torqueConfig}
\end{equation}
It is easy to check that the following matrix
\begin{eqnarray}
\G:=
\left[ \begin{array}{c} \F_a \\ \T_a \end{array} \right]
\label{contrallable}
\end{eqnarray}
is full row rank matrix. As a matter of fact, in engineer 
practice, thruster configuration should always be designed to be 
able to fully control the translation and attitude operations. 
Therefore, we may make the following assumption in the rest of 
the paper:
\vspace{0.05in}
\newline
{\bf Assumption 2:} The configuration matrix $\G$ is always a 
full row rank matrix.
\vspace{0.03in}

\section{Control system design using model predictive control}

Although it is not always straightforward to analyze the close-loop 
stability for MPC control system designs, the following theorem (cf. 
\cite[pages 117-119]{rugh93}) provides a nice sufficient stability 
criterion for the linear time-varying system. 

\begin{theorem}
Suppose for a closed-loop linear time-varying system 
$\dot{\x} = \A(t) \x$ with $\A(t)$ continuously differentiable 
there exist finite positive constants $\alpha$, $\mu$ such that, 
for all $t$, $\| \A(t) \| \le \alpha $ and every pointwise 
eigenvalue of $\A(t)$ satisfies $\Re [\lambda (t)] \le -\mu$.
Then there exists a positive constant $\beta$ such that if the
time derivative of $\A(t)$ satisfies 
$\| \dot{\A}(t) \| \le \beta$ for all $t$, the state equation 
is uniformly exponentially stable.
\label{rughTimeVarying}
\end{theorem}

This theorem is the theoretical base for us to use the MPC
design for the linear time-varying system. 
It is also worthwhile to indicate the following facts: the robust pole 
assignment design can easily achieve the condition that the 
closed-loop system at every fixed time satisfies 
$\Re [\lambda (t)] \le -\mu$, which is
required in Theorem \ref{rughTimeVarying}; for any time between 
the fixed sampling time, robust pole assignment design minimizes 
the sensitivity of the closed-loop poles to the parameter changes
\cite{psnyst14}; moreover, it reduces the impact of disturbance 
torques \cite{yang14}.
Among various pole assignment algorithms, we choose the one 
proposed in \cite{yt93,tity96} because it needs the shortest
computation time among several other popular algorithms 
\cite{psnyst14}, a critical requirement for MPC design. 

We will divide the control force into two parts. The first part
is used to cancel $\n_d(t)$ in (\ref{rendezvousModel2}). This 
can be achieved simply by solving the following linear system of 
equations. 
\begin{equation}
\left[ \begin{array}{c} \F_a \\ \T_a \end{array} \right] \u_1= 
\G \u_1= \left[ \begin{array}{c}
\n_t(t) \\ \n_r(t) \end{array}  \right],
\end{equation}
which gives
\begin{equation}
\u_1= \G^{\dagger} \left[ \begin{array}{c}
\n_t(t) \\ \n_r(t) \end{array}  \right]:= \G^{\dagger} \n,
\label{controlNd}
\end{equation}
where $\G^{\dagger}$ is pseud-inverse of $\G$. In our example,
equations (\ref{forceConfig}) and (\ref{torqueConfig}) implies
$\G^{\dagger} = \G^{-1}$.

The design of second part of the thruster force $\u_2$ is based 
on the following linear time-varying system:
\begin{eqnarray}
\dot{\x} = \A(t) \x + \B \u_2,
\label{rendezvousModel3}
\end{eqnarray}
where $\x$, $\A(t)$, and $\B$ are defined as in 
(\ref{rendezvousModel2}). At every sampling time $t$, $\A(t)$
is evaluated based on the measurable variables. The robust pole 
assignment algorithm of \cite{tity96} is called to get the 
feedback matrix 
\[
\u_2 =  \K(t) \x.
\]
The feedback force $\f_a= \u_1 + \u_2$ is applied to the 
linear time-varying system (\ref{rendezvousModel2}).
The new variables are measured and the next $\A(t)$ is
evaluated in the next sampling time,
and the process is repeated. To avoid the overshoot of relative 
position and relative attitude in the rendezvous and docking 
process, the closed-loop poles should be assigned on the 
negative real axis, i.e., all the poles should be negative and 
real. 

The MPC algorithm using robust pole assignment is summarized as 
follows:

\begin{algorithm} {\ } \\ 
Data: $\mu$, $m_c$, $L_1$, $L_2$, $L_3$, $\J_c$, $\J_t$, $\F_a$, $\T_a$, and $\B$.  
\hspace{0.1in} {\ } \\
Initial condition: At time $t_0$, take the measurements 
${\gamma}={\gamma}_0$, $\r_c$, $\r_t$, $\bar{\q}_{i,tb}$, 
$\bar{\q}_{i,cb}$, $\boldsymbol{\omega}_{i,cb}^{cb}$, 
$\boldsymbol{\omega}_{i,tb}^{tb}$, calculate $\p$, $\r$, 
$\q$, $\R_{tb}^{cb}$, $ \boldsymbol{\omega}$, 
which gives $\x=\x_0$.   {\ } \\
\begin{itemize}
\item[] Step 1: Update $\n_t(r_c,r_t)$, 
$\n_r(\boldsymbol{\omega},\q)$ which gives $\n_d(t)$; update
$\A(t)$ using $\D_t( \dot{\gamma},\ddot{\gamma},\r_c)$, 
$\C_t(\dot{\gamma})$, $\C_r(\boldsymbol{\omega},\q)$, and
$\T(\q)$.
\item[] Step 2: Calculate the gain $\K$ for the linear time-varying
system (\ref{rendezvousModel3}) using robust pole assignment
algorithm implemented as {\tt robpole} (cf. \cite{tity96}).
\item[] Step 3: Apply the controlled thruster force  
$\f_a=\u_1+ \u_2= \G^{\dagger} \n+\K \x$ to 
(\ref{rendezvousModel2}).
\item[] Step 4:  Take the measurements 
${\gamma}$, $\r_c$, $\r_t$, $\bar{\q}_{i,tb}$, 
$\bar{\q}_{i,cb}$, $\boldsymbol{\omega}_{i,cb}^{cb}$, 
$\boldsymbol{\omega}_{i,tb}^{tb}$, calculate $\p$, $\r$, 
$\q$, $\R_{tb}^{cb}$, $ \boldsymbol{\omega}$, 
which gives $\x$. Go back to Step 1.
\end{itemize}
\label{onLine}
\end{algorithm}

The most important step of the algorithm is to call a robust
pole assignment algorithm {\tt robpole} (cf. \cite{tity96,yt93}). 
Let the desired closed-loop eigenvalue set be
$\boldsymbol{\Lambda}=\diag(\lambda_i)$ and 
its corresponding eigenvector matrix be
$\X=[ \x_1, \ldots, \x_n]$ with $\| \x_i \| =1$ such that
\begin{equation}
( \A+\B \K ) \X =  \det(\X) \boldsymbol{\Lambda}.
\label{LambdaX}
\end{equation}
The algorithm of {\tt robpole} is to maximize robustness measurement
of $\det(\X)$ which measures the size of box spanned by the 
eigenvector matrix $\X$ \cite{psnyst14}. That is 
\begin{equation}
\max_{\| \x_i \|=1, i =1,\ldots,n} \X  \hspace{0.1in}
s.t. \hspace{0.1in}  ( \A+\B \K ) \X = \X \boldsymbol{\Lambda}.
\label{optimi}
\end{equation}
The algorithm can be summarized as follows:

\begin{algorithm} {\tt robpole} {\ } \\ 
Data: $\A$, $\B$, and diagonal matrix 
$\boldsymbol{\Lambda}=\diag(\lambda_i)$
with $\lambda_i$ being the desired closed-loop poles.  
\begin{itemize}
\item[] Step 1: QR decomposition for $\B$ yields orthogonal 
$\Q= \left[ \Q_0 \hspace{0.1in} \Q_1 \right]$ 
and triangular $\R$ such that
\begin{equation}
\B= \left[ \Q_0 \hspace{0.1in} \Q_1 \right] \left[ \begin{array} {c}
\R \\ \0  \end{array} \right] .
\label{qrB}
\end{equation}
\item[] Step 2: QR decomposition for $(\A^{\Tr}(t) -\lambda_i \I)\Q_1$
yields orthogonal $\V= \left[ \V_{0i} \hspace{0.1in} \V_{1i} \right]$ 
and triangular $\Y$ such that
\begin{equation}
(\A^{\Tr}(t) -\lambda_i \I)\Q_1= \left[ \V_{0i} 
\hspace{0.1in} \V_{1i} \right] 
\left[ \begin{array} {c}
\Y \\ \0  \end{array} \right] ,
\hspace{0.1in} i=1, \ldots,n.
\label{qrV}
\end{equation}
\item[] Step 3: Cyclically select one real or a pair of (real or 
complex conjugate) unit length eigenvectors such that 
$\x_i \in \mathcal{S}_i = span(\V_{1i})$ and the robustness measure
$\det(\X)$ is maximized.
\item[] Step 4: The feedback matrix is given by
\begin{equation}
\K =\R^{-1}\Q_0^{\Tr}(\X \boldsymbol{\Lambda} \X^{-1} - \A(t)).
\label{robK}
\end{equation}
\end{itemize}
\label{onLine1}
\end{algorithm}

\begin{remark}
It is worthwhile to emphasize that $\B$ in 
(\ref{rendezvousModel3}) is a constant matrix.
Therefore, Step 1 can be calculated off-line, which saves
a lot of computational burden for the MPC control scheme.
\end{remark}

In \cite{yang14}, it is shown that maximizing $\det(\X)$
is equivalent to minimizing an upper-bound of $\H_{\infty}$
norm of the closed-loop system. Therefore, the design 
is insensitive to modeling error. Moreover, maximizing  $\det(\X)$
amounts to minimizing the an upper-bound of the gain of the channel 
from disturbance input to the output. This can be seen as follows:
Assume that system (\ref{rendezvousModel3})
with disturbance $\f_d$ is 
\[
\dot{\x} = \A(t) \x + \B \u + \f_d,
\]
and its output vector is defined by 
\[
\y = \x..
\]
Since $\u = \K \x$, taking the Laplace transformation gives
\[
\Y(s) = (s\I-(\A+\B\K))^{-1} \f_d(s) = \X (s\I -\Lambda )^{-1}
\X^{-1} \f_d(s).
\] 
This gives
\begin{equation}
\| \Y(s) \| \le  \| (s\I-\Lambda )^{-1} \| \cdot 
\| \X \| \cdot \| \X^{-1} \|  \cdot  \|  \f_d(s) \|
\label{bound}
\end{equation}
where $\| (s\I-\Lambda )^{-1} \|$ is a constant,
$ \|  \f_d(s) \|$ is the unknown disturbance, and 
$\kappa :=\| \X \| \cdot \| \X^{-1} \| =\sigma_1/\sigma_n$
($\sigma_1$ and $\sigma_n$ are the maximum and minimum
singular values of $\X$). According to \cite{yang14}, 
maximizing $\det(\X)$ implies 
minimizing the upper bound of $\kappa=\| \X \| \cdot \| \X^{-1} \|$.
Therefore, from (\ref{bound}), we conclude that
the robust pole assignment design reduces the
impact of $\f_d$ on the output of $\Y$, which means that
the design has good disturbance rejection property.

Applying Theorem \ref{rughTimeVarying}, Yang \cite{yang18} provides
some moderate conditions for uniformly exponential stability for a
class of the MPC designs, which can easily be extended to Algorithm \ref{onLine}. 

\begin{theorem}[\cite{yang18}]
Let the prescribed closed-loop eigenvalues be distinct and
$\Re [\lambda_i (t)] \le -\mu$ for some $\mu >0$, 
and the time dependent variables in LTV system $(\A(t),\B(t))$
and their derivatives be bounded. Then, the closed-loop 
time-varying system (\ref{rendezvousModel}) with the robust 
pole assignment control obtained by Algorithm \ref{onLine}
is uniformly exponentially stable.
\label{MPCmain}
\end{theorem}

\section{Simulation test}

In this section, we present an simulation test example to support
the design idea. We use the simulation example of \cite{zd12} and
all the parameters used in that paper. We compare the simulation
results of the two designs to demonstrate the superiority of the
proposed design.

First, the physics constants, such as, gravitational constant
$\mu=3.986004418*10^{14} m^3 /(kg\cdot s^2)$, Earth radius 
$6371000 \mbox{ m}$, are from \cite{wertz78}. The rest parameters
are from \cite{zd12}: the target spacecraft orbit is circular and
the altitude is $250$ km, $L_1 =L_2 = L_3=2$ m, the mass of the 
chaser is $10$ kg and its inertia matrix is 
$\J_c=\diag [10,\,\, 10,\,\, 10] kg \cdot m^2$, the mass of the 
target $10$ kg and its inertia matrix is given as
\[
\J_t=\left[ \begin{array}{ccc}
10   &   2.5  & 3.5 \\
2.5  &   10   & 4.5 \\
3.5  &   4.5  & 10
\end{array} \right] kg \cdot m^2,
\] 
$\F_a$ is given in (\ref{forceConfig}), $\T_a$ is given in 
(\ref{torqueConfig}). The initial condition is set as
$\p(0)=[10,\,\, -10,\,\, 10]^{\Tr} m$, 
$\d(0)=[5,\,\,  -4,\,\, 4]^{\Tr} m/s$,
$\bar{\q}(0)=[0.3772,\,\, -0.4329,\,\, 0.6645,\,\, 0.4783]^{\Tr}$,
$\boldsymbol{\omega}(0)=[0,\,\, 0,\,\, 0]^{\Tr} rad/s$. 

\begin{figure}[hb!]
\centering
\includegraphics[height=8cm,width=8cm]{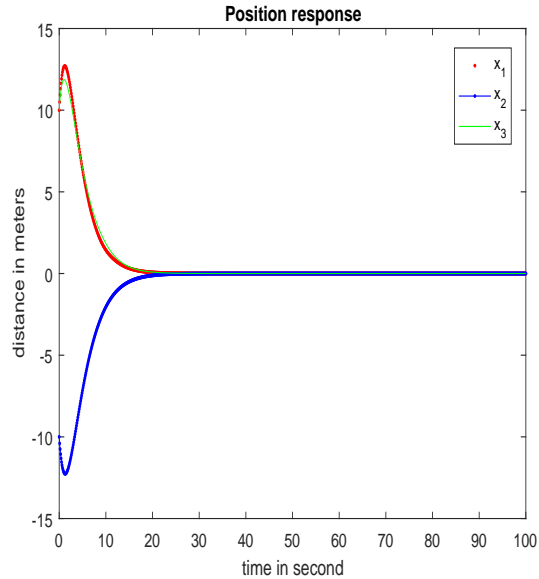}
\caption{Position response.}
\label{fig:PositionResponse1}
\end{figure}

\begin{figure}[ht]
\centering
\includegraphics[height=8cm,width=8cm]{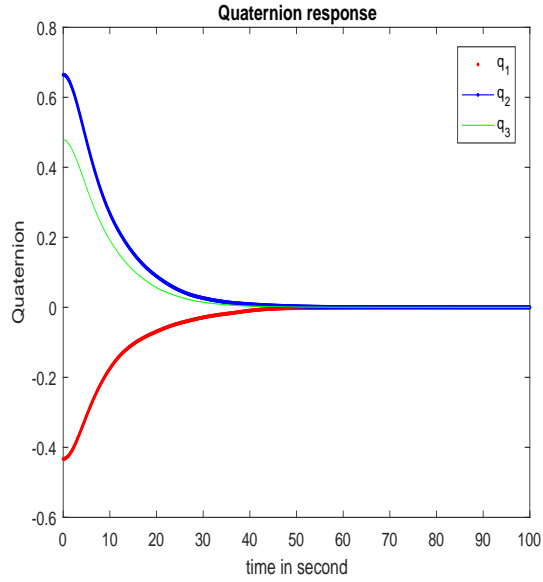}
\caption{Attitude response.}
\label{fig:AttitudeResponse1}
\end{figure}

\begin{figure}[hb!]
\centering
\includegraphics[height=8cm,width=8cm]{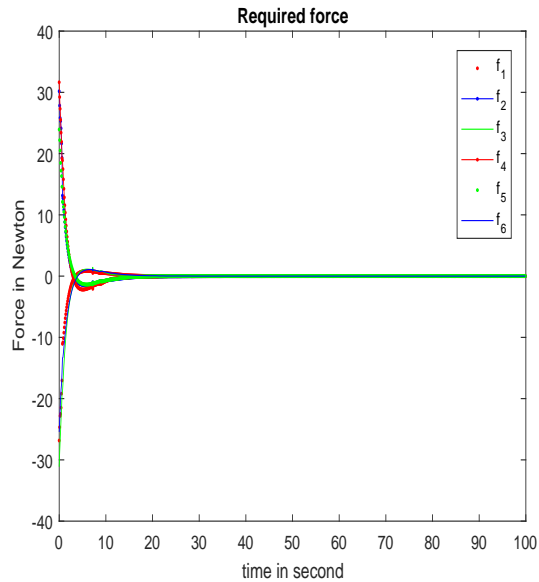}
\caption{Required forces.}
\label{fig:forceRequired1}
\end{figure}

To avoid the overshoot in relative distance and relative
attitude to guarantee the soft docking, we should avoid the
oscillation. Therefore, the proposed closed loop poles are
set to
\[
-0.1 -0.2 -0.15 -0.25 -0.3 -0.35 -0.4 -0.45 -0.5 -0.55 -0.6 -0.65
\]
Applying the on-line algorithm \ref{onLine} to this problem again,
we performed the simulation of the closed-loop response. The 
results are shown in Figures \ref{fig:PositionResponse1} -
\ref{fig:forceRequired1}. Fig. \ref{fig:PositionResponse1} is the 
response of relative position between the chaser and the target 
and Fig. \ref{fig:AttitudeResponse1} is the response of relative 
attitude between the chaser and the target. These figures show 
that the design successfully avoid the oscillation in the 
docking process and achieved the soft docking.
Fig. \ref{fig:forceRequired1} is the forces in $6$ thrusters 
used in this docking process, the maximum forces is about
$31$ Newton, which 
is smaller than the maximum forces used in the design of 
\cite{zd12}, which is in the range of $360$ Newton.

Comparing to the simulation tests in \cite{ssj06,gys09,zld11,zwld14,
sh15,zd12,kng08,slf11,xjlky15}, the simulation using the proposed
method is the only one that does not have oscillation in all these responses, 
which is a clear indication that the design achieves soft docking.

\begin{figure}[ht]
\centering
\includegraphics[height=8cm,width=8cm]{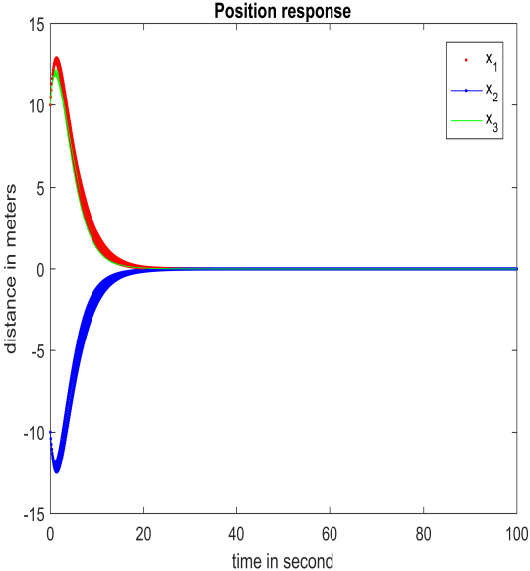}
\caption{Position response for randomly generated 100 problems.}
\label{fig:position100Runs}
\end{figure}

\begin{figure}[hb!]
\centering
\includegraphics[height=8cm,width=8cm]{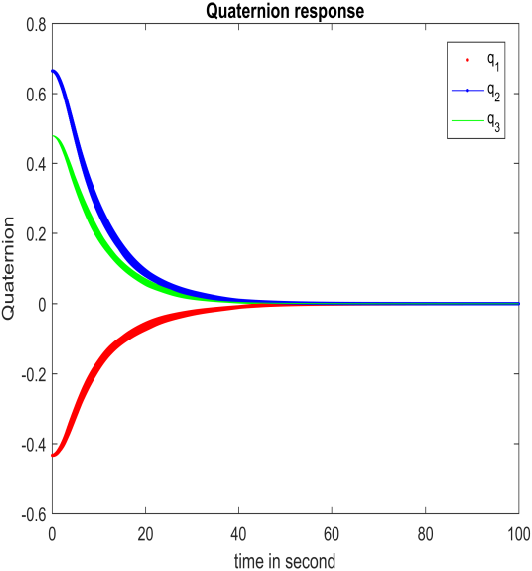}
\caption{Attitude response for randomly generated 100 problems.}
\label{fig:attitude100Runs}
\end{figure}

To verify that the design is insensitive to the modeling error on 
the performance and the design reduces the impact of disturbance
torques, we consider the first set of prescribed closed-loop 
eigenvalues with (a) randomly generated modeling errors in 
all entries of $\J_t$ and $\J_c$ for up to $1kg \cdot m^2$,
(b) randomly generated installation errors in thruster location
$L_1$, $L_2$, and $L_3$ for up to $0.01m$, and (c) randomly
generated disturbance with mangetude for up to $10 \%$
$\B * \F * \| \x(t) \| $ for all $t>0$. The system is controlled
by using the controller designed in this section without the
information on the randomly generated modeling and installation
error and randomly generated disturbances. The performance of 
the position responses and atttitude responses for $100$
randomly generated problems described above is provided
in Figures \ref{fig:position100Runs} and \ref{fig:attitude100Runs}.
Clearly, the performance of position and attitude responses is 
insensitive to the modeling error and disturbance torques.

\section{Conclusions} 
 
In this paper, we proposed an MPS-based design method for 
translation and attitude combined control in spacecraft
rendezvous and soft docking. A robust pole assignment design
is used in the MPS-based method, which enjoys several nice
features: (a) all closed-loop poles are negative real which 
guarantees the strict requirement of soft docking, (b) on-line
computation is moderate and all intermediate solutions are
feasible, (c) the design is insensitive to the effects of the modeling 
error on the performance, and (d) the design reduces the impact
of disturbance torques. A simulation example is included to 
demonstrate the effectiveness and the efficiency of the proposed
design.


\end{document}